\documentclass{amsart}

\usepackage{graphicx}
\usepackage{latexsym, amsmath, amsfonts, amscd, amssymb, verbatim}
\usepackage{enumerate}
\usepackage{enumitem}
\usepackage{cancel}
\usepackage{graphicx}
\usepackage{float}
\usepackage{tikz,pgfplots}

\usetikzlibrary{calc} 
\usetikzlibrary{arrows.meta,bending,decorations.markings,intersections} 
\tikzset{
	curve arrow/.style args={%
		to pos #1 with length #2}{
		decoration={
			markings,
			mark=at position 0 with {\pgfextra{%
					\pgfmathsetmacro{\tmpArrowTime}{#2/(\pgfdecoratedpathlength)}
					\xdef\tmpArrowTime{\tmpArrowTime}}},
			mark=at position {#1-\tmpArrowTime} with {\coordinate(@1);},
			mark=at position {#1-2*\tmpArrowTime/3} with {\coordinate(@2);},
			mark=at position {#1-\tmpArrowTime/3} with {\coordinate(@3);},
			mark=at position {#1} with {\coordinate(@4);
				\draw[-{Stealth[length=#2,bend]}]       
				(@1) .. controls (@2) and (@3) .. (@4);},
		},
		postaction=decorate,
	}
}

\newtheorem{theorem}{Theorem}[section]
\newtheorem{lemma}[theorem]{Lemma}
\newtheorem{corollary}[theorem]{Corollary}
\newtheorem{proposition}[theorem]{Proposition}

\theoremstyle{definition}
\newtheorem{definition}[theorem]{Definition}
\newtheorem{example}[theorem]{Example}

\theoremstyle{remark}
\newtheorem{remark}[theorem]{Remark}

\numberwithin{equation}{section}

\begin{document}

\title{A PROOF OF Lee-Lee's CONJECTURE ABOUT GEOMETRY OF RIGID MODULES}

\author{Son Dang Nguyen}
\address{Department of Mathematics, The University of Alabama, Tuscaloosa, USA }
\email{sdnguyen1@crimson.ua.edu}

\subjclass[2000]{ 13F60 (16G20)}

\date{\today}


\keywords{rigid modules, real Schur roots, cluster algebras, quiver representations}

\begin{abstract}
This paper proves Lee-Lee's conjecture that establishes a coincidence between the set of associated roots of non-self-intersecting curves in a $n$-punctured disc and the set of real Schur roots of acyclic (valued) quivers with $n$ vertices. 
\end{abstract}

\maketitle

\section{Introduction}
Given an algebraically closed field $\mathbb{F}$ and an acyclic quiver $Q$ with $n$ vertices, an indecomposable module $M$ is called by \textbf{rigid} if $Ext_{\mathbb{F}Q}^1(M,M) = 0$ over the path algebra $\mathbb{F}Q$. It is well-known that from the quiver $Q$ one may construct the irreducible \textbf{symmetric} Cartan matrix companion $A=A(Q)$, the root system $\Phi  = \Phi (A)$ and the Coxeter-Weyl groups $G_n = G_n(A)$ with ordered simple reflections $\left\{ {{s _1},...,{s_n}} \right\}$.  In the root systems,  positive \textbf{real Schur roots} are dimension vectors of indecomposable rigid modules. Such rigid modules and their real Schur roots play very important roles in representation theory and cluster algebras, see \cite{AFPT}, \cite{ST}, \cite{Kac1}, \cite{Kac2},  \cite{KIRS}, \cite{SZ} just to name, hence it leads to necessarily characterize them among all real roots of $\Phi$. Remark that the real Schur roots depend on the orientation of the quiver $Q$, which makes it quite difficult to characterize them.  In \cite{LL}, K. H. Lee and K. Lee suggested a conjecture about their geometric property.
\begin{flushleft}
	\textbf{Conjecture}(\cite{LL}). \textit{The real roots assigned by simple curves are precisely real Schur roots. (defined in Section 2). }	
\end{flushleft}
In the paper, they proved it for acyclic quivers of $3$ vertices with multiple arrows between every pair of vertices. When the paper was published on arxiv, A. Felikson and P. Tumarkin \cite{AFPT} proved it for all acyclic quivers of finite vertices with multiple arrows between every pair of vertices (called by $2$\textit{-complete} acyclic quivers). In the paper, we show the solution of the conjecture for all acyclic (valued) quivers of finite vertices.

Before going to the proof, we importantly need to reform the conjecture in another version. Given a quiver $Q$, we have a bijection between the orientation of the quiver and an admissible-sink (or admissible-source) order of its Coxeter element $c:=s_{i_1}s_{i_2}...s_{i_n}$. Hence without any risk, in the paper we can fix the quiver $Q$ such that its Coxeter element is $c=s_1s_2...s_n$. In \cite{KIRS}, Igusa and Schiffler showed that reflections corresponding real Schur roots are precisely prefix reflections in the factorization of the Coxeter element $c$. In \cite{AH}, Hubery and Krause gave a characterization of real Schur roots in terms of simple partition. Their reflections are reflection elements belonging to the poset of generalized simple partitions $NC=NC(Q):=\left\{ {w \in G_n\left| {1 \le w \le c} \right.} \right\}$ where $\le$ denotes the \textit{absolute order} on $G_n$. 

In more general cases for a valued quiver $Q$, given $A$ be an irreducible \textbf{symmetrizable} generalized Cartan matrix with its symmetrizer $D$ and orientation $O$, in \cite{GLS1} Geiß, Leclerc and Schröer define a Iwanaga-Gorenstein $\mathbb{F}$-algebra $H:=H_\mathbb{F}(A,D,O)$ for any field $\mathbb{F}$ in terms of a quiver with relations with a hereditary algebra $\tilde H$ of the corresponding type and a Noetherian $\mathbb{F[[\epsilon]]}$-algebra ${\hat H}$. In particular, in \cite{GLS2}, their main results showed that the indecomposable \textit{rigid locally free} $H$-modules are parametrized via their rank vector, by the real Schur roots associated to $(A,O)$. Moreover, the \textit{left finite bricks} of $H$ are parametrized via their dimension vector by the real Schur roots associated to $(C^T,O)$. Also as the symmetric case of acyclic quivers $Q$, real Schur roots were proven to be real roots in its associated root system such that their reflections are reflection elements of the simple partition $NC$ in its associated Weyl group. From the observations,  Lee-Lee's conjecture may be reformed as follow:
\begin{flushleft}\textbf{Conjecture}(\cite{LL})
	Three statements are equivalent:
	\begin{enumerate}
		\item A positive real root $\beta$ is Schur.
		\item Its reflection $r_\beta$ is a prefix of the Coxeter element $c$ i.e. there exists $n-1$ reflections $r_2,r_3,...,r_n$ such that $r_\beta r_2,...,r_n=c$.
		\item The real root $\beta$  may be presented by a simple curve  and its  reflection may be presented by a simple closed curve.
	\end{enumerate}
\end{flushleft}
where \cite{KIRS} and \cite{GLS2} proved $(1) \Leftrightarrow (2)$, so we only need to show $(2) \Leftrightarrow (3)$. 

Section $2$ is devoted to recall the reader the constructions used before giving the proof of Lee-Lee's conjecture in section $3$. In section $4$, we present another proof on finite types that also implies an one-side proof for affine types. We also give some other results and open problems that may be completed in the future.
 
\begin{flushleft}
	\textbf{Acknowledgements:} This work is a part of my PhD thesis, supervised by Proofessor Kygungong Lee. I would like to thank him for his guidance and patience. The author is also deeply grateful to The University of Alabama for providing ideal working conditions to write this note.
\end{flushleft}
\section{Reminders}

Given an irreducible symmetrizable Cartan matrix $A=(a_{ij}) \in M_{{\rm{n}} \times {\rm{n}}}(\mathbb{Z})$, it is well-known that one may construct its associated root system $\Phi=\Phi(A)$ as follows. We fix its simple roots $\{\alpha_1,...,\alpha_n\}$ and define its simple reflections $\{s_1,...,s_n\}$  by 
\[{s_i}({\alpha _j}) = {\alpha _j} - {a_{ij}}{\alpha _i},\begin{array}{*{20}{c}}
	{}&{i,j = 1,...,n}
\end{array}.\]
These reflections generates a group ${G_n}={G_n}(A) = \left\langle {{s_1},...,{s_n}\left| {s_i^2 = 1} \right.,{{({s_i}{s_j})}^{{m_{ij}}}} = 1} \right\rangle$ called by the Coxeter$/$Weyl group where $m_{ij}$ are defined from the table: 
\begin{center}
	\begin{tabular}{ |c|c|c|c|c|c| } 
		\hline
		$a_{ij}a_{ji}$ & 0 & 1&2&3&$\ge$ 4 \\
		\hline
		$m_{ij}$ & 2 & 3&4&6&$\infty$ \\ 
		\hline
	\end{tabular}
\end{center}
Let $B_n$ be the braid group on $n$ strands and abstractly presented
\[{B_n} = \left\langle {{\sigma _1},...,{\sigma _{n - 1}}\left| {{\sigma _i}{\sigma _{i + 1}}{\sigma _i} = {\sigma _{i + 1}}{\sigma _i}{\sigma _{i + 1}},{\sigma _i}{\sigma _j} = {\sigma _j}{\sigma _i}} \right.\text{ if }\left| {i - j} \right| > 1} \right\rangle \]
with the action on $n$ copies of an arbitrary group $H$ as follow:
\[\begin{array}{l}
{\sigma _i}.({g_1},...,{g_i},{g_{i + 1}},...,{g_n}) = ({g_1},...,{g_i}{g_{i + 1}}{g_i}^{-1},{g_i},...,{g_n}),\\
\sigma _i^{ - 1}.({g_1},...,{g_i},{g_{i + 1}},...,{g_n}) = ({g_1},...,{g_{i + 1}},{g_{i + 1}}^{-1}{g_i}{g_{i + 1}},...,{g_n})
\end{array}\]
for any sequence $({g_1},...,{g_n}) \in {H^n}$. It is clear that the action fixes the product $g_1g_2...g_n$. This implies that the position of reflections appearing in the factorization of the Coxeter element of $G_n$ may be ignored by the following lemma. 
\begin{lemma}
	Assume that a reflection $r_\beta$ belonging a factorization of the Coxeter element $c$, then it is a prefix of $c$ i.e. there exists $n-1$ reflections $r_2,r_3,...,r_n$ such that $r_\beta r_2,...,r_n=c$. 
\end{lemma}
\textbf{Proof.} The smallest number of factorization of $c$ into reflections is $n$ (see \cite{KIRS}). Assume that $c=r_1r_2...r_n$ and $r_j=r_\beta$ for some $2 \le j\le n$, then $r_\beta$ is a prefix of $\sigma_2^{-1}...\sigma_{j}^{-1}.(r_1,r_2,...,r_n)$. The action of $B_n$ keeps the product invariant, hence the proof is completed.  $\blacksquare$

Now we will present real roots as curves and their reflections as closed curves of by the following construction of D. Bessis in \cite{DB}. For convenience we will draw the \textit{disc} $D \subset \mathbb{C}$ as an upper half-plane $\{ z \in {\Bbb C}\left| {\operatorname{Im} z > 0} \right.\} $ with the \textit{punctured point} $p_i$ placed from left to right on the horizontal line $Imz=1$ and the point $\infty$ at infinity. Denote by ${\ell _i}$ a vertical ray ${\ell _i} = \left\{ {z\left| {\operatorname{Re} z = \operatorname{Re} {p_i},\operatorname{Im} z > \operatorname{ Im} {p_i}} \right.} \right\}$. 

\begin{definition}
	A curve $\gamma$ is a continuous map $\gamma:[0,1] \to D$ such that $\gamma(0)=O,\gamma(1) \in \{p_1,p_2,...,p_n\}$ and $\gamma(t) \notin \{O,p_1,p_2,...,p_n\}$ for $t\in(0,1)$. 
	
	A \textit{simple curve} is a non-self-intersecting curve. 
	
A closed curve $\hat \lambda$ is a continuous map $\hat \lambda:[0,1] \to D$ such that $\gamma(0)=\gamma(1)=O$ and $\gamma(t) \notin \{O,p_1,p_2,...,p_n\}$ for $t\in(0,1)$. 

A non-self-intersecting closed curve is called by a \textit{simple loop}.

Denote $Ins(\hat \lambda)$ be the interior region bounded by the closed curve $\hat \lambda$. 
\end{definition}

\begin{remark}
	One may associate such a curve $\gamma$ uniquely  (up to isotopy) to a closed curve $\hat{\gamma}$ as follows: starting from $\gamma(0)=O$, follow $\gamma$; arriving close to $\gamma(1)$, make a possible turn around a small circle centred on $\gamma(1)$; return to $O$ following $\gamma$ backwards. Thus a simple curve associates to a simple loop and so the set of simple curves may inject as a subset of the set of simple loops.
\end{remark}

\begin{definition}
Given a closed curve $\hat\lambda$, then the word $w_{\hat\lambda}$ presents an element of $G_n$ constructed as follows: following $\hat\lambda$, write $s_j$ each time it crosses some rays ${\ell _j}$.

\begin{definition}
	Given a curve $\gamma$ with $\gamma(1)=p_i$ and its closed curve $\hat \gamma$, then the word $\alpha_\gamma$ presents a real root constructed as follows: following $\gamma$, write $s_j$ each time it crosses some rays ${\ell _j}$ and add $\alpha_i$ as last word. It is clear that the word $w_{\hat \gamma}$ presents its reflection, denoted by $r_{\hat\gamma}$.
\end{definition}

\end{definition}
From the definition of curves and closed curves, now we may identify them and their words that present real roots and elements in $G_n$. Since all real roots and elements in $G_n$ may be presented by words, the set of real roots is the set of curves (up to isotopy) and $G_n$ is the set of closed curves (up to isotopy). Note that two curves may be isotopic together but the former presents a positive real root $\alpha$, the latter presents its negative real root $-\alpha$ as the following example.
\begin{example}
	
\begin{figure}[h!]
	\begin{tikzpicture}[line width=1pt]
	\node  at (0.9,5) [label=below:{$\ell_1$}]{};
	\node  at (2,5) [label=below:{$\ell_2$}]{};		
	\node  at (3,5) [label=below:{$\ell_3$}]{};	
	
	\node at (0,0) [circle,fill,inner sep=0pt,minimum size=3pt,label=below:{$O$}]{};
	\node  at (0.9,2) [circle,fill,inner sep=0pt,minimum size=3pt,label=below:{$p_1$}]{};
	\node  at (2,2) [circle,fill,inner sep=0pt,minimum size=3pt,label=below:{$p_2$}]{};		
	\node  at (3,2) [circle,fill,inner sep=0pt,minimum size=3pt,label=below:{$p_3$}]{};
	
	\node  at (-0.5,1) [label=below:{$\hat \lambda$}]{};
	\node  at (1,0.5) [label=below:{$\gamma$}]{};	
	\node  at (1.2,0) [label=below:{$\hat \gamma$}]{};		
	
	\draw plot [smooth] coordinates {
		(0,0) (1,0.5) (1.7,2.8) (3,2)   
	};
	\draw [red] plot [smooth] coordinates {
		(0,0) (-0.5,2) (0.5,3.5) (4,3.5) (3,-1) (0,0)
	};
	\draw [green] plot [smooth] coordinates {
		(0,0) (0.8,0.5) (1.5,3) (3.2,2.3) (3.2,1.8) (1.8,2.4) (1.2,0) (0,0)
	};
	
	\end{tikzpicture}
\end{figure}
In the example of $G_3$, the simple red loop $\hat \lambda$ induces the word $s_1s_2s_3$ of the Coxeter element $c$. In particular, the loop is isotopic to the boundary $\partial D = {S^1}$ of the disc $D$. The simple dark curve $\gamma$ presents a positive real root $\alpha_\gamma=s_2\alpha_3$ and its simple green loop $\hat{\gamma}$ presents a reflection $r_{\hat\gamma}=s_2s_3s_2$. From the curve $\gamma$ we may draw another simple curve $\beta$ that is isotopic with $\gamma$ by going around $p_3$ one time before ending at $p_3$. Then $\beta$ presents a negative real root $\alpha_\beta=s_2s_3\alpha_3=-\alpha_\gamma$. 
\end{example}

A real root of $\Phi$ and an element of $G_n$ may be presented by many different words i.e. curves and closed curves up to isotopy. However, if ${G_n} = {W_n}: = \{ {s_1},...,{s_n}\left| {s_1^2 = ... = s_n^2 = 1} \right.\} $, called by \textit{universal Coxeter groups}, then  real roots of $\Phi$ and elements of $G_n$ are determined uniquely by their words because Cayley graphs of $W_n$ are trees. Hence real roots and elements of $G_n$ are presented uniquely by curves and closed curves. Now we prove main results that lead to the Lee-Lee's conjecture. 
\begin{lemma}
	Given two simple curves with distinct end points such that they do not intersect each other except for $O$. We order them in a clock wise order of their emanating from $O$, written by the sequence $({\gamma _1},{\gamma _2})$. Then
	\begin{enumerate}
		\item The product ${r_{{\hat\gamma _1}}}{r_{{\hat\gamma _2}}}$ may be presented by a simple loop $\hat\gamma$ such that two simple loops ${{{\hat\gamma _1}}},{{{\hat\gamma _2}}} \subset Ins(\hat\gamma)$.
		\item Braid group $B_2$ acts on $({r_{{\hat\gamma _1}}},{r_{{\hat\gamma _2}}})$ to be ${\sigma _1}.({r_{{\hat\gamma _1}}},{r_{{\hat\gamma _2}}}) = ({r_{{\hat\gamma _1}}}{r_{{\hat\gamma _2}}}{r_{{\hat\gamma _1}}},{r_{{\hat\gamma _1}}})$ (similarly for ${\sigma _1}^{-1}.({r_{{\hat\gamma _1}}},{r_{{\hat\gamma _2}}}) = ({r_{{\hat\gamma _2}}}, {r_{{\hat\gamma _2}}}{r_{{\hat\gamma _1}}}{r_{{\hat\gamma _2}}})$), then it is possible to construct a simple curve $\gamma$ presenting the real root ${r_{{\hat\gamma _1}}}{\alpha_{{\gamma _2}}}$ (thus the simple loop $\hat\gamma$ presenting ${r_{{\hat\gamma _1}}}{r_{{\hat\gamma _2}}}{r_{{\hat\gamma _1}}}$) such that $\gamma$ does not intersect with $\gamma_2$ (thus $\hat\gamma$ does not intersect with $\hat\gamma_2$) and $\gamma(1)=\gamma_2(1)$. Hence braid group action preserves non-intersecting property of simple curves and loops. 
		  
	\end{enumerate}	
\end{lemma} 
\textbf{Proof.} For $(1)$, since two simple curves $({{{\gamma _1}}},{{{\gamma _2}}})$ do not intersect each other, neither do their simple loops $({{\hat\gamma _1},\hat{\gamma _2}}) $. Therefore we can connect them in a small neighbor at $O$ to become the loop $\hat\gamma$ qualified as the follow example.

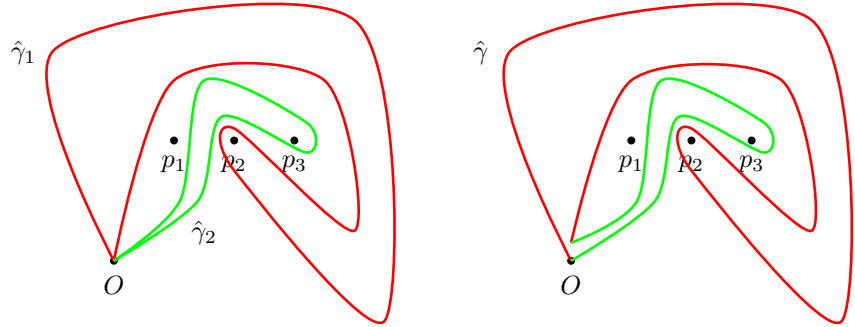
\begin{figure}[h!]
	\begin{center}
	\begin{tikzpicture}[scale=0.8,line width=1pt]
	\node at (0,0) [circle,fill,inner sep=0pt,minimum size=3pt,label=below:{$O$}]{};
	\node  at (1,2) [circle,fill,inner sep=0pt,minimum size=3pt,label=below:{$p_1$}]{};
	\node  at (2,2) [circle,fill,inner sep=0pt,minimum size=3pt,label=below:{$p_2$}]{};		
	\node  at (3,2) [circle,fill,inner sep=0pt,minimum size=3pt,label=below:{$p_3$}]{};
	
	\node  at (1.5,1) [label=below:{$\hat \gamma_2$}]{};	
	\node  at (-1.5,4) [label=below:{$\hat \gamma_1$}]{};	
	
	\draw [red] plot [smooth] coordinates {
		(0,0) (-1,3.5) (4,4) (4.5,-1) (2,1.5) (2,2.2) (4,0.5) (3.5,3) (1,3) (0,0)
	};

	\draw [green] plot [smooth] coordinates {
	(0,0) (1.1,1) (1.5,3) (3.2,2.3) (3.2,1.8) (1.8,2.4) (1.4,1) (0,0)
};
	\end{tikzpicture}
\qquad 
	\begin{tikzpicture}[scale=0.8,line width=1pt]
	\node at (0,0) [circle,fill,inner sep=0pt,minimum size=3pt,label=below:{$O$}]{};
	\node  at (1,2) [circle,fill,inner sep=0pt,minimum size=3pt,label=below:{$p_1$}]{};
	\node  at (2,2) [circle,fill,inner sep=0pt,minimum size=3pt,label=below:{$p_2$}]{};		
	\node  at (3,2) [circle,fill,inner sep=0pt,minimum size=3pt,label=below:{$p_3$}]{};
	
	\node  at (-1.5,4) [label=below:{$\hat \gamma$}]{};	
	
	\draw [red] plot [smooth] coordinates {
		(0,0) (-1,3.5) (4,4) (4.5,-1) (2,1.5) (2,2.2) (4,0.5) (3.5,3) (1,3) (0,0.3)
	};
	\draw [green] plot [smooth] coordinates {
		(0,0.3) (1.1,1) (1.5,3) (3.2,2.3) (3.2,1.8) (1.8,2.4) (1.4,1) (0,0)
	};
	
	\end{tikzpicture}
		\end{center}
\end{figure}

\newpage
For $(2)$, we present $({r_{{\hat\gamma _1}}}{r_{{\hat\gamma _2}}}{r_{{\hat\gamma _1}}},{r_{{\hat\gamma _2}}})$ from $({r_{{\hat\gamma _1}}},{r_{{\hat\gamma _2}}})$ as follows: in a small neighbor of O we connect the loop $\hat \gamma _1$ with the curve $\gamma _2$, then we get a new simple curve $\gamma'$ presenting the real roots ${r_{{\hat\gamma _1}}}{{\alpha_{\gamma _2}}}$ as qualified and its corresponding simple loop $\hat\gamma'$ presents ${r_{{\hat\gamma _1}}}{r_{{\hat\gamma _2}}}{r_{{\hat\gamma _1}}}$. Finally, we reorder two curve $\gamma_1$ and $\gamma'$ in clock wise order of emanating from O, then we have the sequence $(\hat\gamma',\hat\gamma_1)$ presenting for $({r_{{\hat\gamma _1}}}{r_{{\hat\gamma _2}}}{r_{{\hat\gamma _1}}},{r_{{\hat\gamma _1}}})$. The remaining part $({r_{{\hat\gamma _2}}}, {r_{{\hat\gamma _2}}}{r_{{\hat\gamma _1}}}{r_{{\hat\gamma _2}}})$ is similar.  $\blacksquare$

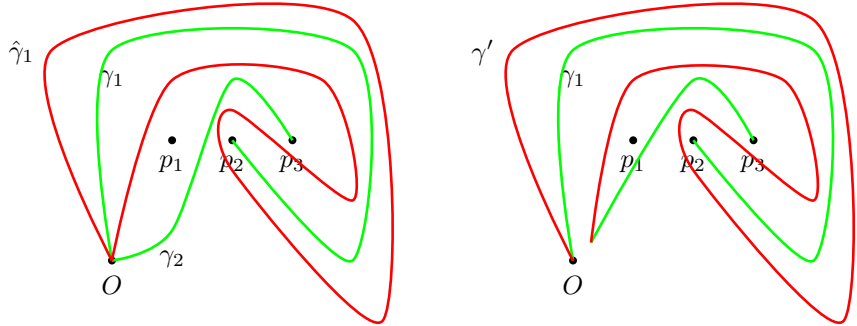
	\begin{figure}[h!]
		\begin{center}	
	\begin{tikzpicture}[scale=0.8,line width=1pt]
	\node at (0,0) [circle,fill,inner sep=0pt,minimum size=3pt,label=below:{$O$}]{};
	\node  at (1,2) [circle,fill,inner sep=0pt,minimum size=3pt,label=below:{$p_1$}]{};
	\node  at (2,2) [circle,fill,inner sep=0pt,minimum size=3pt,label=below:{$p_2$}]{};		
	\node  at (3,2) [circle,fill,inner sep=0pt,minimum size=3pt,label=below:{$p_3$}]{};
	
	\node  at (1,0.5) [label=below:{$\gamma_2$}]{};	
	\node  at (0,3.5) [label=below:{$\gamma_1$}]{};	
	\node  at (-1.5,4) [label=below:{$\hat \gamma_1$}]{};	
			
	\draw [green] plot [smooth] coordinates {
		(0,0) (1,0.5) (2,3) (3,2)   
	};
	\draw [green] plot [smooth] coordinates {
	(0,0) (0,3.5) (4,3.5) (4,0) (2,2)
};
	\draw [red] plot [smooth] coordinates {
		(0,0) (-1,3.5) (4,4) (4.5,-1) (2,1.5) (2,2.5) (4,1) (3.5,3) (1,3) (0,0)
	};
			 
	\end{tikzpicture}
	\qquad
	\begin{tikzpicture}[scale=0.8,line width=1pt]
	\node at (0,0) [circle,fill,inner sep=0pt,minimum size=3pt,label=below:{$O$}]{};
	\node  at (1,2) [circle,fill,inner sep=0pt,minimum size=3pt,label=below:{$p_1$}]{};
	\node  at (2,2) [circle,fill,inner sep=0pt,minimum size=3pt,label=below:{$p_2$}]{};		
	\node  at (3,2) [circle,fill,inner sep=0pt,minimum size=3pt,label=below:{$p_3$}]{};
	
	\node  at (0,3.5) [label=below:{$\gamma_1$}]{};	
	\node  at (-1.5,4) [label=below:{$\gamma'$}]{};
	
	\draw [green] plot [smooth] coordinates {
		(0.3,0.3) (2,3) (3,2)   
	};
	\draw [green] plot [smooth] coordinates {
		(0,0) (0,3.5) (4,3.5) (4,0) (2,2)
	};
	\draw [red] plot [smooth] coordinates {
		(0,0) (-1,3.5) (4,4) (4.5,-1) (2,1.5) (2,2.5) (4,1) (3.5,3) (1,3) (0.3,0.3)
	};
	
	\end{tikzpicture}
\end{center}
\end{figure}

\begin{corollary}
	Given $n$ simple curves with distinct end points such that their pairs do not intersect each other except for $O$. We order them in a clock wise order of emanating from $O$, possibly written by the sequence $({\gamma _1},...,{\gamma _n})$. Then
	\begin{enumerate}
		\item All products ${r_{{\hat\gamma _i}}}...{r_{{\hat\gamma _j}}},1\le i \le j \le n$ may be presented by loops $\hat\gamma_{ij}$ such that  $j-i+1$ simple loops ${{{\hat\gamma _i}}},...,{{{\hat\gamma _j}}} \subset Ins(\hat\gamma_{ij})$ and $Ins(\hat\gamma_{ij}) \subset Ins(\hat\gamma_{im})$ for $1 \le i \le j \le m \le n$. 
		
		In particular, the loop $\hat\gamma_{1n}$ presents the product ${r_{{\hat\gamma _1}}}...{r_{{\hat\gamma _n}}}$, all simple loops ${{{\hat\gamma _1}}},...,{{{\hat\gamma _n}}} \subset Ins(\hat\gamma_{1n})$ and $Ins(\hat\gamma_{ij}) \subset Ins(\hat\gamma_{1n})$ for $1 \le i \le j \le n$.
		\item Braid group $B_n$ acting on $({r_{{\hat\gamma _1}}},...,{r_{{\hat\gamma _n}}})$ preserves non-intersecting property of simple curves and loops.
	\end{enumerate}
\end{corollary}	
\textbf{Proof.} For $(1)$, the proof is similar with $(1)$ of Lemma 2.7 by applying it for $({\gamma _i},{\gamma _{i+1}}),1\le i\le j-1\le n$. For $(2)$, we only need to check actions of ${\sigma _i}$, $1\le i\le n-1$
on $({r_{{\hat\gamma _1}}},...,{r_{{\hat\gamma _n}}})$. But the actions only locally transform two loops in the sequence and fix the remaining $n-2$ loops, so its proof is completed from $(2)$ of Lemma 2.7. $\blacksquare$
\newpage
\begin{proposition}
	The simple loop $\hat\gamma_{1n}$ is isotopic with the simple loop presenting the Coxeter element $c=s_1...s_n$. Hence ${r_{{\hat\gamma _1}}}...{r_{{\hat\gamma _n}}}=c$.
\end{proposition}
\textbf{Proof.} The idea of the proof is the same as Lemma 2.3 and 2.4 in \cite{DB}. Corollary 2.8 implies that $Ins(\hat\gamma_{1n})$ contains all punctured points $\{p_1,...,p_n\}$ and presents their product ${r_{{\hat\gamma _1}}}....{r_{{\hat\gamma _n}}}$. In the other hand, $c$ may be presented by a loop $\hat\gamma_c$ with the exact word that contains $\{p_1,...,p_n\}$ and it is isotopic with $\partial D = {S^1}$. Thus it may be chosen in its isotopic class such that  $Ins(\hat\gamma_{1n}) \subset Ins(c)$. Since the annulus between $\hat\gamma_c$ and $\hat\gamma_{in}$ contains no punched point in $\{p_1,...,p_n\}$, they are isotopic and so ${r_{{\hat\gamma _1}}}...{r_{{\hat\gamma _n}}}=c$. $\blacksquare$ 

\section{Proof of Lee-Lee's conjecture}
\begin{theorem}
	Three statements are equivalent:
\begin{enumerate}
	\item A real root $\beta$ is Schur.
	\item Its reflection $r_\beta$ is a prefix of the Coxeter element $c$ i.e. there exists $n-1$ reflections $r_2,r_3,...,r_n$ such that $r_\beta r_2,...,r_n=c$.
	\item The real root $\beta$  may be presented by a simple curve  and its  reflection may be presented by a simple closed curve.
\end{enumerate}
where \cite{KIRS} and \cite{GLS2} proved $(1) \Leftrightarrow (2)$, so we only need to prove $(2) \Leftrightarrow (3)$. 
\end{theorem}
\textbf{Proof.} A real Schur root $\alpha$ may be presented by a simple curve. Indeed, in \cite{KIRS} and \cite{GLS2}, a real root $\alpha$ is Schur if and only if its reflection $r_\alpha$ is a prefix of the Coxeter element $c$. Therefore it may be induced in a reflection sequence $(r_1,...,r_n)$ where $r_\alpha=r_j$ with some $j$ such that $r_1...r_n=c$. Since the braid group $B_n$ acts transitively on the factorization of $c$ (see \cite{KIRS}), there exists $\sigma \in B_n$ such that $\sigma.(s_1,...,s_n)=(r_1,...,r_n)$. It is clear that $(s_1,...,s_n)$ can be presented by $n$ simple loops without pairwise intersection each other. But Corollary 2.8 shows that the action of $B_n$ preserves non-intersecting property of $n$ simple loops presenting them, so the real Schur root $\alpha$ may be presented by a simple curve.

Conversely, a simple curve can be induced in $n$ simple curves with no pairwise intersection because of induction on $n$ by cutting the disc $D$ along the curve giving rise to an $n-1$-punctured disc. We order them in a clock wise order of emanating from O, then Proposition 2.9 shows that the product of their reflections is the Coxeter element $c$. Thus roots corresponding reflections of these simple curve are Schur roots because of \cite{KIRS} and \cite{GLS2}. $\blacksquare$

\begin{corollary}
	In the case $G_n=W_n:=\left\langle {{s_1},...,{s_n}\left| {s_1^2 = ... = s_n^2} \right.} \right\rangle $, we have a bijective correspondence between simple curves (up to isotopy) and positive real Schur roots. 
\end{corollary}
\textbf{Proof.} In the case, real roots and reflections are uniquely presented by their reduced words, thus non-isotopic curves present distinct real roots. $\blacksquare$

This is the same result obtained in \cite{AFPT}. 

Now we let $\Phi ^{NC}$ be the set of simple loops and define an equivalence $\sim$ as follows: $\hat \gamma    \sim \hat \lambda$ if they present the same elements in $G_n$ and denote by ${{\tilde \Phi }^{NC}}:=\Phi ^{NC}/\sim$.

Remark that any two isotopic loops present the same elements but conversely, it is not true because a reflection may be presented by many non-isotopic loops. 

\begin{definition}
	For $\forall \hat \gamma ,\hat \lambda  \in {{\tilde \Phi }^{NC}}$ a partial order $\subseteq $ in ${{\tilde \Phi }^{NC}}$ is defined by
	$$\hat \gamma  \subseteq \hat \lambda \mathop  \Leftrightarrow \limits^{def}  \hat \gamma, \hat \lambda \text{ may be chosen in their equivalent class such that } Ins(\hat \gamma) \subseteq Ins(\hat \lambda)$$ and number of punctured points in $Ins(\hat \gamma)$ is fewer than in $ Ins(\hat \lambda)$.
\end{definition}

The \textit{absolute length} $l(w)$ of $w\in G_n$ is the minimal $k \ge 0$ such that $w$ can be written as product $w=t_1t_2...t_r$ of reflections $t_i\in G_n$.

\begin{definition}
	For $\forall w,u \in G_n$ an absolute order $\le$ on $NC$ is defined by 
	$$w \le u  \mathop  \Leftrightarrow \limits^{def} l(w)+l(w^{-1}u)=l(u).$$
\end{definition}

Recall the simple partition $NC:=\left\{ {w \in G_n\left| {1 \le w \le c} \right.} \right\}$. 
\begin{theorem}
	We have an order-preserving isomorphism between $({{\tilde \Phi }^{NC}},\subseteq)$ and $(NC,\le)$.
\end{theorem}
\textbf{Proof.} Bijection between ${{\tilde \Phi }^{NC}}$ and $NC$ is trivial from definition of ${{\tilde \Phi }^{NC}}$. Given $ \hat \gamma ,\hat \lambda  \in {{\tilde \Phi }^{NC}}$ such that $\hat \gamma  \subseteq \hat \lambda {\text{ }}$ and $w_{\hat \gamma}$, $w_{\hat \lambda} \in G_n$ are elements that they present for, respectively. We may assume that all punctured points in $Ins(\hat \lambda)$ are $\{p_{i_1},...,p_{i_k},p_{i_{k+1}},...,p_{i_m}\}$ and all punctured points in $Ins(\hat \gamma)$ are $\{p_{i_1},...,p_{i_k}\}$ for $1 \le i_1 \le i_{k} \le i_m \le n$. In the annulus between $\hat \gamma$ and $\hat \lambda$ we draw $m-k$ simple loops $\hat \beta_{i_{j}}$ such that each loop contains exact one punctured point $p_{i_j}$ (thus each one presents a reflection $r_{i_{j}}$) with $k+1 \le j \le m$, they are pairwise non-intersecting, and they do not intersect with $\hat \lambda$. We order $m-k+1$ loops $\{\hat \gamma,\beta_{i_{k+1}},...,\beta_{i_{m}}\}$ in a clock wise order of their eliminating from $O$, then the product of their corresponding elements is $w_{\hat \lambda}$. Hence $w_{\hat \gamma} \le w_{\hat \lambda}$. Similarly, $w_{\hat \lambda} \le c$. Conversely, given $1 \le u \le w \le c$, then non-crossing partition $NC$ implies that there are $n$ reflections $r_{i_{1}},...,r_{i_{n}}$ such that $r_{i_{j}}r_{i_{j+1}}...r_{i_{k}}=u$, $r_{i_{j}}r_{i_{j+1}}...r_{i_{k}}r_{i_{k+1}}...r_{i_{m}}=w$ and $r_{i_{1}}...r_{i_{n}}=c$ for $1 \le j \le k \le m \le n $. From Corollary 2.8, two simple loops $ \hat \gamma ,\hat \lambda  \in {{\tilde \Phi }^{NC}}$ presenting $u,w$ may be chosen such that $\hat \gamma  \subseteq \hat \lambda {\text{ }}.$ $\blacksquare$

Denote $s(c)$ and $t(c)$ be first and last simple reflections of $c$, respectively. Let $s(c)(Q)$ be a new quiver obtained by conversing arrows adjacent to the vertex corresponding to $s(c)$ (similarly for $t(c)(Q)$). The new quivers obtained by this approach correspond to \textbf{mutation} of quivers (see \cite{AFPT}, \cite{NS2}) at sink-source vertices corresponding to $s(c)$ and $t(c)$ in the theory of quiver representations. These mutations maintain root systems of the quivers but change their orientation, thus change their set of real Schur roots. However, the set of real Schur roots of the new quivers and their corresponding simple curves may be obtained by the following proposition. 
\begin{proposition}
	A real root $\beta$ is Schur in the quiver $Q$ if and only if the real root $s(c)\beta$ is Schur in the quiver $s(c)(Q)$ (similarly for $t(c)(Q)$). 
	
	 The word $w_\beta$ is obtained from a simple curve $\gamma_\beta$ in the setting of the quiver $Q$. Then a simple curve presenting the real Schur root $s(c)\beta$ in the setting of the quiver $s(c)(Q)$ may be constructed from the word $s(c)w_\beta$ (similarly for $t(c)(Q)$). 
\end{proposition}
\textbf{Proof.} Since $\beta$ is a real Schur root, there exists $n-1$ reflections $r_1,...,r_{n-1}$ such that $r_\beta r_1...r_{n-1}=c$, thus $(s(c)r_\beta s(c))(s(c)r_1s(c))...s(c))(s(c)r_{n-1}s(c))=(s(c)cs(c))$. The right hand side is the Coxeter element corresponding of the new quiver $s(c))(sQ)$ so the first part of the proposition is proved. For the latter part, we may consider it when the Weyl group of the quiver $Q$ is an universal Coxeter group. In the case, all reduced words obtained from simple curves are unique, thus simple curve presenting for $s(c)\beta$ has to be constructed from the word $s(c)\beta$. $\blacksquare$

\begin{example}	
We give an example for a rank-$3$ quiver  with its Coxeter element $c:=s_1s_2s_3$, so $s(c)=s_1$ and $t(c)=s_3$. A real Schur root $\beta:=s_2s_3\alpha_2$ is presented by an simple curve $\gamma_\beta$. Mutating the quiver at the vertex $1$ corresponding to $s(c)$ we obtain the new quiver $s(c)(Q)$ with its Coxeter element $\bar c:=s_2s_3s_1=s(c)cs(c)$.
 
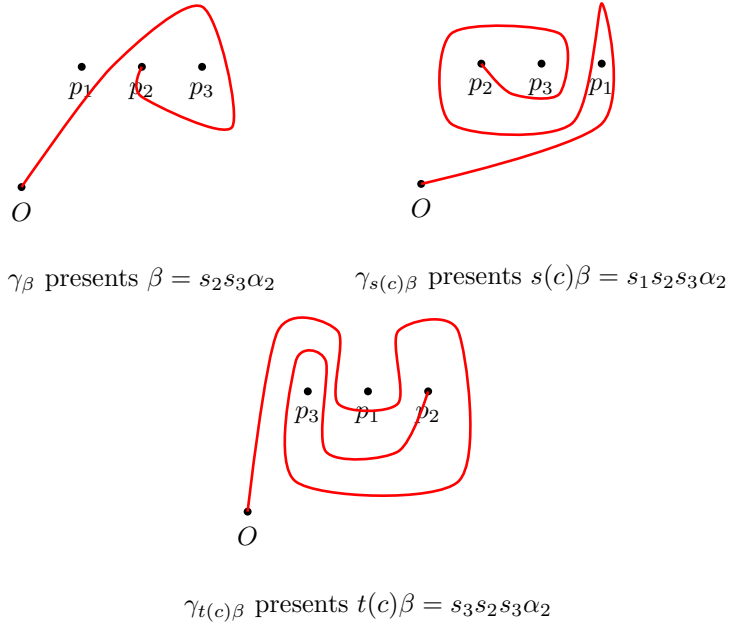
\begin{figure}[h!]
	\begin{center}
		\begin{tikzpicture}[scale=0.8,line width=1pt]
			\node at (0,0) [circle,fill,inner sep=0pt,minimum size=3pt,label=below:{$O$}]{};
			\node  at (1,2) [circle,fill,inner sep=0pt,minimum size=3pt,label=below:{$p_1$}]{};
			\node  at (2,2) [circle,fill,inner sep=0pt,minimum size=3pt,label=below:{$p_2$}]{};		
			\node  at (3,2) [circle,fill,inner sep=0pt,minimum size=3pt,label=below:{$p_3$}]{};

			\node  at (2,-1) [label=below:{$\gamma_\beta$ presents $\beta=s_2s_3\alpha_2$}]{};	
			
			\draw [red] plot [smooth] coordinates {
				(0,0) (1.5,2) (3,3) (3.5,1) (2,1.5) (2,2) 
			};		
		\end{tikzpicture}
		\qquad 
		\begin{tikzpicture}[scale=0.8,line width=1pt]
			\node at (0,0) [circle,fill,inner sep=0pt,minimum size=3pt,label=below:{$O$}]{};
			\node  at (1,2) [circle,fill,inner sep=0pt,minimum size=3pt,label=below:{$p_2$}]{};
			\node  at (2,2) [circle,fill,inner sep=0pt,minimum size=3pt,label=below:{$p_3$}]{};		
			\node  at (3,2) [circle,fill,inner sep=0pt,minimum size=3pt,label=below:{$p_1$}]{};

			\node  at (2,-1) [label=below:{$\gamma_{s(c)\beta}$ presents $s(c)\beta=s_1s_2s_3\alpha_2$}]{};

			\draw [red] plot [smooth] coordinates {
			(0,0) (3,1) (3,3) (2.5,1) (0.5,1) (0.5,2.5) (2.3,2.5) (2.3,1.5) (1.5,1.5) (1,2)  
			};
			
		\end{tikzpicture}
			\qquad 
	\begin{tikzpicture}[scale=0.8,line width=1pt]
		\node at (0,0) [circle,fill,inner sep=0pt,minimum size=3pt,label=below:{$O$}]{};
		\node  at (1,2) [circle,fill,inner sep=0pt,minimum size=3pt,label=below:{$p_3$}]{};
		\node  at (2,2) [circle,fill,inner sep=0pt,minimum size=3pt,label=below:{$p_1$}]{};		
		\node  at (3,2) [circle,fill,inner sep=0pt,minimum size=3pt,label=below:{$p_2$}]{};
		
		\node  at (2,-1) [label=below:{$\gamma_{t(c)\beta}$ presents $t(c)\beta=s_3s_2s_3\alpha_2$}]{};

		\draw [red] plot [smooth] coordinates {
			(0,0) (0.5,3) (1.5,3) (1.5,1.8) (2.5,1.8) (2.5,3) (3.5,3) (3.5,0.5) (0.8,0.5) (0.8,2.5) (1.3,2.5) (1.3,1) (2.5,1) (3,2)
		};
		
	\end{tikzpicture}
	\end{center}
\end{figure}
\end{example}

\section{Some Remarks on Finite, Affine and rank-$2$ Types}
In the section, we give another proof for the Lee-Lee's conjecture in the case of finite and affine types. It also yields an algorithm to construct simple curves for all real Schur roots in finite, affine and rank-$2$ types. Assume $G_n$ is a Weyl group of finite types with the Dynkin diagram corresponding $\Delta$. 

\begin{lemma} \label{coxeteraction}
	Given a simple curve $\gamma$ with its real root $\alpha_\gamma$ and its reflection $r_\gamma$, then the real roots $c\alpha_\gamma$ and $c^{-1}\alpha_\gamma$ may be presented by simple curves with the exact words. Hence the action of the Coxeter element $c$ preserves non-self-intersecting property of simple curves. This yields the equivalence between $c$-orbit of simple curves and  $c$-orbit of their real roots.	
\end{lemma}
\textbf{Proof.} The proof is straight-forward from spiraling the simple curves  clockwise or counterclockwise as the figures:   
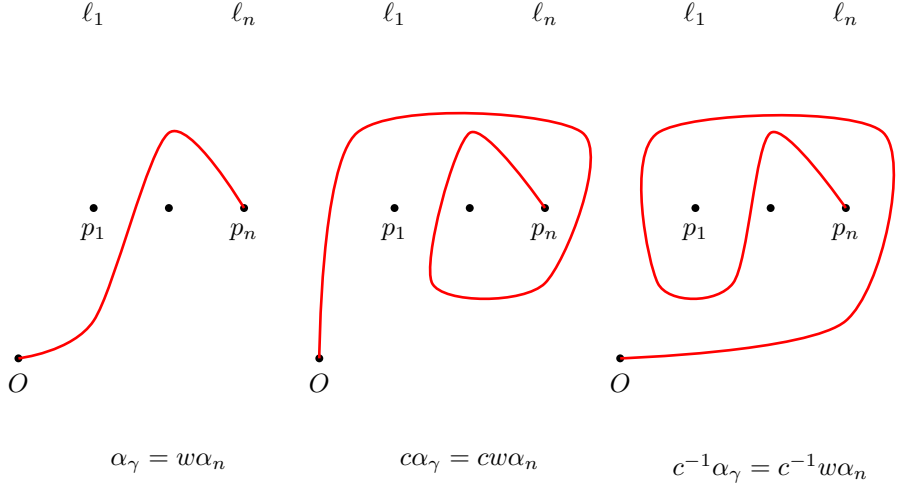
\begin{figure}[h!]
	\begin{tikzpicture}[line width=1pt]
	\node at (0,0) [circle,fill,inner sep=0pt,minimum size=3pt,label=below:{$O$}]{};
	\node  at (1,2) [circle,fill,inner sep=0pt,minimum size=3pt,label=below:{$p_1$}]{};
	\node  at (1,5) [label=below:{$\ell_1$}]{};	
	\node  at (2,2) [circle,fill,inner sep=0pt,minimum size=3pt,label=below:{}]{};		
	\node  at (3,5) [label=below:{$\ell_n$}]{};	
	\node  at (3,2) [circle,fill,inner sep=0pt,minimum size=3pt,label=below:{$p_n$}]{};
	\node  at (2,-1) [label=below:{$\alpha_\gamma  = w\alpha_n$}]{};			
	\draw [red] plot [smooth] coordinates {
		(0,0) (1,0.5) (2,3) (3,2)   
	};
	
	\node at (4,0) [circle,fill,inner sep=0pt,minimum size=3pt,label=below:{$O$}]{};
	\node  at (5,2) [circle,fill,inner sep=0pt,minimum size=3pt,label=below:{$p_1$}]{};
	\node  at (5,5) [label=below:{$\ell_1$}]{};		
	\node  at (6,2) [circle,fill,inner sep=0pt,minimum size=3pt,label=below:{}]{};	
	\node  at (7,5) [label=below:{$\ell_n$}]{};	
	\node  at (7,2) [circle,fill,inner sep=0pt,minimum size=3pt,label=below:{$p_n$}]{};
	\node  at (6,-1) [label=below:{$c\alpha_\gamma  = cw\alpha_n$}]{};
	\draw [red] plot [smooth] coordinates {
		(4,0) (4.5,3) (7.5,3) (7,1) (5.5,1) (6,3) (7,2)
	};
	
	\node at (8,0) [circle,fill,inner sep=0pt,minimum size=3pt,label=below:{$O$}]{};
	\node  at (9,2) [circle,fill,inner sep=0pt,minimum size=3pt,label=below:{$p_1$}]{};
	\node  at (9,5) [label=below:{$\ell_1$}]{};		
	\node  at (10,2) [circle,fill,inner sep=0pt,minimum size=3pt,label=below:{}]{};	
	\node  at (11,5) [label=below:{$\ell_n$}]{};	
	\node  at (11,2) [circle,fill,inner sep=0pt,minimum size=3pt,label=below:{$p_n$}]{};
	\node  at (10,-1) [label=below:{$c^{-1}\alpha_\gamma  = c^{-1}w\alpha_n$}]{};
	\draw [red] plot [smooth] coordinates {
		(8,0) (11,0.5) (11.5,3) (8.5,3) (8.5,1) (9.5,1) (10,3) (11,2)
	};
	
	\end{tikzpicture}
\end{figure}

\newpage
\begin{remark}
	Since $	{({\sigma _{n-1}}{\sigma _{n-2}}...{\sigma _{1}})^n}.({s_1},...,{s_n}) = (c{s_1}{c^{ - 1}},...,c{s_n}{c^{ - 1}})$ and \\
	${({\sigma _{n-1}}{\sigma _{n-2}}...{\sigma _{1}})^{-n}}.({s_1},...,{s_n}) = (c^{ - 1}{s_1}{c},...,c^{ - 1}{s_n}{c})$, the action of $c$ on a simple loop may be seen as a special case of the action of the braid group $B_n$.
	
	Lemma 4.1 also shows a connection of \textit{Auslander-Reitein translation} $\tau$ in acyclic quiver representation theory and the Coxeter element $c$. While the Auslander-Reitein translation preserves the rigid property of a rigid module (see \cite{IDA}), the Coxeter element preserves non-self-intersecting property of their corresponding real roots.
\end{remark}
 
\begin{corollary}
	For $n=2$, Lee-Lee's conjecture holds and the set of real Schur roots is $c$-orbit of $\left\{ {{\alpha _1},{\alpha _2}} \right\}$ which contains all real roots. Therefore in the case all positive real roots are Schur.
\end{corollary}
\textbf{Proof.} 
 From \cite{KIRS}, \cite{GLS2} and the transitive property of the action of $B_2$ on $(s_1,s_2)$, one knew that reflections of real Schur roots are precisely elements appear in the orbit of $B_2$ on $(s_1,s_2)$. In particular, for $B_2=<\sigma_1>$, we have
\[\begin{array}{l}
	\sigma _1^{2h}({s_1},{s_2}) = ({c^h}{s_1}{c^{ - h}},{c^h}{s_2}{c^{ - h}}) \text{ for }h \in \mathbb{Z},\\
	\\
	\sigma _1^{2k + 1}({s_1},{s_2}) = ({c^k}{s_1}{s_2}{s_1}{c^{ - k}},{c^k}{s_1}{c^{ - k}})\text{ for }k \in \mathbb{Z^+},\\
	\\
	\sigma _1^{ - (2k + 1)}({s_1},{s_2}) = ({c^{ - k}}{s_2}{c^k},{c^{ - k}}{s_2}{s_1}{s_2}{c^k}), \text{ for }k \in \mathbb{Z^+}.
\end{array}\]
Moreover, $s_1s_2s_1=cs_2c^{-1}$ and $s_2s_1s_2=c^{-1}s_1c$. Hence Lemma 4.1 completed the proof. $\blacksquare$

\begin{corollary}
	Lee-Lee's conjecture holds for finite-type root systems.
\end{corollary}
\textbf{Proof.} In finite-type cases, real roots are precisely real Schur roots, so the proof of simple curves presenting real Schur roots is trivial. Conversely, let ${\beta _k} = {s_1}{s_2}...{s_{k - 1}}{\alpha _k}$, for $1\le k \le n$, then it is clear that these roots may be presented by simple curves with the exact word. Moreover, finite-type root systems have exact $n$ distinct $c$-orbits where ${\beta _k} $ belongs to the each one (see Proposition 33, chapter VI in \cite{N}). Thus Lemma 4.1 implies Lee-Lee's conjecture. $\blacksquare$ 

In the case of affine types, the action of $c$ also gives an one-side proof of Lee-Lee's conjecture and the remain ones is delivered from Section 3.
\begin{corollary}
	Lee-Lee's conjecture holds for affine-type root systems.
\end{corollary}
\textbf{Proof.} In \cite{NS1} and \cite{NS2}, the authors show explicit $\Phi _c^{re}$ set of real Schur roots as follows: the set has exact  $2n$ infinite $c$-orbits and $n-2$ finite $c$-orbits with $n \ge 2$. The transversal set of $2n$ infinite $c$-orbits includes ${\beta _k} = {s_1}{s_2}...{s_{k - 1}}{\alpha _k}$ and ${\delta _k} = {s_n}{s_{n - 1}}...{s_{k + 1}}{\alpha _k}$, for $1\le k \le n$ that clearly might be presented by simple curves. Moreover, all real roots in $n-2$ finite $c$-orbits are of finite types that might be presented by simple curves because of Corollary 4.2. Hence Lemma 4.1 implies that all real Schur roots of affine types may be presented by simple curves. $\blacksquare$  

Let $R$ be the set of reflections of $G_n$ and recall the action of $B_n$ on $s=(s_1,...,s_n) \in R^n$. We denote $c=c(\Delta)$ and $h=h(\Delta)$ respectively being the Coxeter element and the Coxeter number of Dynkin diagram $\Delta$. Let $B_{G_n}$ be the stabilizer subgroup of $s$ and $B_n(s)$ is the $B_n$-orbit at $s$. Transitivity action of $B_n$ on the completed exceptional sequence implies that the orbit is all possibly completed exceptional sequences.

\begin{proposition}
 	We have a bijection between the coset $B_n/B_{G_n}$ and the factorization of the Coxeter element $c$, hence the index $[B_n:B_{G_n}]=\frac{{n!h{^n}}}{{\left| {{G_n}} \right|}}$.

\end{proposition}
\textbf{Proof.} Bijection is trivial from the orbit-stabilizer theorem and the specific formula of number of the factorization of the Coxeter group is shown in \cite{ONSRM}. $\blacksquare$

The following table in \cite{ONSRM} exhibits the index formula for the connected Dynkin diagrams $\Delta$:

\begin{center}
\begin{tabular}{l|llllllll}
$\Delta$	& $\mathbb{A}_n$ & $\mathbb{B}_n,\mathbb{C}_n$ & $\mathbb{D}_n$ & $\mathbb{E}_6$  & $\mathbb{E}_7$  & $\mathbb{E}_8$  & $\mathbb{F}_4$  &$\mathbb{G}_2$  \\ 
	\hline
$\frac{{n!h{^n}}}{{\left| {{G_n}} \right|}}$&$(n+1)^{n-1}$  &$n^n$  &$2(n-1)^n$  & $2^9.3^4$ & $2.3^{12}$  &$2.3^5.5^7$  &$2^4.3^3$  &$2.3$  \\ 
\end{tabular}
\end{center}
Note that $\left\{ {{\sigma_1},{\sigma_1}{\sigma_2}...{\sigma_{n-1}}} \right\} $ generates $B_n$ and the group $B_{G_n}$ is a finitely generated subgroup because it is a finite index subgroup of $B_n$. Unfortunately, now we cannot yet find the generating set of $B_{G_n}$ but introduce some elements of $B_{G_n}$. Let $\Delta(i--j)$ to be subdiagrams of the Dynkin diagram $\Delta$ with vertices $\{i,i+1,...,j\}$ and $h_{ij}:=h(\Delta(i--j))$ to be their corresponding Coxeter number; particularly $h_i:=h_{ii+1}$.
\begin{lemma}
	For $1\le i < j \le n-1$, 
		$${{{({\sigma _{j-1}}{\sigma _{j-2}}...{\sigma _{i}})}^{(j-i+1){h_{ij}}}},{{({\sigma _{n-1}}...{\sigma _{1}})}^{nh}} \in {B_{G_n}}}.$$
		In particular, if $(s_is_{i+1})^{m_{ij}}=1$, then $\sigma_i^{m_{ij}} \in B_{G_n}$ with $m_{i,j}=2,3,4,6$.
\end{lemma}
\textbf{Proof.} We have 
	\[\begin{array}{*{20}{l}}
	\begin{array}{l}
		\sigma _i^2.({s_1},...{s_i},{s_{i + 1}},...{s_n}) = ({s_1},...,{c_i}{s_i}{c_i}^{ - 1},{c_i}{s_{i + 1}}{c_i}^{ - 1},...{s_n}),\\
		{({\sigma _{j-1}}{\sigma _{j-2}}...{\sigma _{i}})^(j-i+1)}.({s_1},...{s_i},...,{s_{i + 1}},...,{s_n}) = ({s_1},...{c_{ij}}{s_i}c_{ij}^{ - 1},...,{c_{ij}}{s_{i + 1}}c_{ij}^{ - 1},...,{s_n}),
	\end{array}\\
	{({\sigma _{n-1}}{\sigma _{n-2}}...{\sigma _{1}})^n}.({s_1},...,{s_n}) = (c{s_1}{c^{ - 1}},...,c{s_n}{c^{ - 1}}),
\end{array}\]
where $c_{i}=c(\Delta(i--i+1))=s_is_{i+1}$, $c_{ij}=c(\Delta(i--j))=s_is_{i+1}...s_{j}$. Since $c_i^{h_i}=c_{ij}^{h_{ij}}=c^h=id$, the proof is completed. $\blacksquare$

Remark that ${({\sigma _{n-1}}{\sigma _2}...{\sigma _{1}})^{nh}}$ belongs to the center of $B_n$, hence ${({\sigma _{n-1}}{\sigma _2}...{\sigma _{1}})^{nh}} \in N_n:=\bigcap\limits_{\sigma  \in {A_n}} {\sigma {B_{G_n}}{\sigma ^{ - 1}}} $ i.e. $B_n$ acts unfaithfully on $s$. Since $B_{G_n}$ is a finite-index subgroup, so is $N_n$.
\begin{example}
	For $G_2=\mathbb{A}_2$ we have $[B_2:B_{\mathbb{A}_2}]=3$ and $B_{\mathbb{A}_2}=<{\sigma _1}^{3}>$.
	
	For $G_2=\mathbb{G}_2$ we have $[B_2:B_{\mathbb{G}_2}]=6$ and $B_{\mathbb{G}_2}=<{\sigma _1}^{6}>$.
\end{example}

\begin{corollary}
	Assume that there exists a proper normal subgroup $M$ that is strictlly larger than $B_{G_n}$ for $n \ge3$, then $n=2$ $(mod$ $3)$. Hence $B_n$ and ${\raise0.7ex\hbox{${{B_n}}$} \!\mathord{\left/
			{\vphantom {{{B_n}} M}}\right.\kern-\nulldelimiterspace}
		\!\lower0.7ex\hbox{$N_n$}}$ is not solvable for $n \ge 3$. 
\end{corollary}
\textbf{Proof.} We consider the natural projection $\pi :{B_n} \to {\raise0.7ex\hbox{${{B_n}}$} \!\mathord{\left/
		{\vphantom {{{B_n}} M}}\right.\kern-\nulldelimiterspace}
	\!\lower0.7ex\hbox{$M$}}$. Since $\sigma_{1}^3 \in B_{A_n}$ and $M$ is proper, the order of $\pi(\sigma_{1})$ is $3$ in ${\raise0.7ex\hbox{${{B_n}}$} \!\mathord{\left/
		{\vphantom {{{B_n}} M}}\right.\kern-\nulldelimiterspace}
	\!\lower0.7ex\hbox{$M$}}$ so $3$ divides $[B_n:B_{\mathbb{A}_n}]=(n+1)^{n-1}$. Thus $n=2$ $(mod$ $3)$. This implies that there is no proper normal subgroup that is stictlly larger than $B_{\mathbb{A}_3}$ in $B_3$. Moreover, $B_{\mathbb{A}_3}$ is not a normal subgroup of $B_3$,  hence $B_3$ is not solvable. Since $B_3$ may be embedded in $B_n$ for $n \ge 3$, $B_n$ is not solvable and so is ${\raise0.7ex\hbox{${{B_n}}$} \!\mathord{\left/
		{\vphantom {{{B_n}} M}}\right.\kern-\nulldelimiterspace}
	\!\lower0.7ex\hbox{$N_n$}}$. $\blacksquare$

Finally, we finish the paper with several further questions which might be of interest.

\textbf{Question 1:} How can we find the finite generating set of $B_{G_n}$ and $N_n:=\bigcap\limits_{\sigma  \in {A_n}} {\sigma {B_{G_n}}{\sigma ^{ - 1}}} $ when $G_n$ is a Weyl group of finite types?

\textbf{Question 2:} What is the classification of the finite non-solvable group ${\raise0.7ex\hbox{${{B_n}}$} \!\mathord{\left/
		{\vphantom {{{B_n}}}}\right.\kern-\nulldelimiterspace}
	\!\lower0.7ex\hbox{$ N_n$}}$ for $n \ge 3$?

\bibliographystyle{amsplain}

\begin{thebibliography}{10}
	
\bibitem {AH} Andrew Hubery and Henning Krause. \textit{A categorification of simple partitions}. J. Eur. Math. Soc. (JEMS) 18 (2016), no. 10, 2273–2313

\bibitem {AFPT}  Anna Felikson and Pavel Tumarkin. \textit{Acyclic cluster algebras, reflection groups, and curves on a punctured disc}. Adv. Math. 340 (2018), 855–882.

\bibitem{GLS1} Christof Geiß; Bernard Leclerc and Jan Schröer. \textit{Quivers with relations for symmetrizable Cartan matrices I: Foundations}. Invent. Math. 209 (2017), no. 1, 61–158.

\bibitem{GLS2} Christof Geiß, Bernard Leclerc and Jan Schröer. \textit{Rigid modules and Schur roots}. Mathematische Zeitschrift (2019). 

\bibitem {DB} David Bessis. \textit{A dual braid monoid for the free group}. J. Algebra 302 (2006), no. 1, 55–69.

\bibitem {ST} David Speyer and Hugh Thomas, \textit{Acyclic cluster algebras revisited. Algebras, quivers and representations}, 275–298, Abel Symp., 8, Springer, Heidelberg, 2013.

\bibitem {IDA} Ibrahim Assem, Daniel Simson, Andrzej Skowroński, \textit{Elements of the representation theory of associative algebras}. Vol. 1. Techniques of representation theory. London Mathematical Society Student Texts, 65. Cambridge University Press, Cambridge, 2006. x+458 pp. ISBN: 978-0-521-58423-4; 978-0-521-58631-3; 0-521-58631-3

\bibitem {Kac1} V. G. Kac, \textit{Infinite root systems, representations of graphs and invariant theory}. Invent. Math. 56 (1980), 57–92.
 
\bibitem {Kac2} V. G. Kac, \textit{Infinite root systems, representations of graphs and invariant theory. II}. J. Algebra 78 (1982), 141–162.

\bibitem {KIRS}  Kiyoshi Igusa and Ralf Schiffler. \textit{ Exceptional sequences and clusters}. J. Algebra 323 (2010), no. 8, 2183–2202. 

\bibitem {LL} Kyu-Hwan Lee and Kyungyong Lee. \textit{A correspondence between rigid modules over path algebras and simple curves on riemann surfaces}, Exp. Math, 2019.

\bibitem {ONSRM} Mustafa Obaid, Khalid Nauman, Wafa S. M. Al-Shammakh, Wafaa Fakieh, and Claus Michael Ringel, \textit{The number of complete exceptional sequences for a Dynkin algebra}. Colloq. Math. 133 (2013), no. 2, 197–210.

\bibitem {NS1} Nathan Reading and Salvatore Stella.\textit{ The action of a Coxeter element on an affine root system}. Proc. Amer. Math. Soc. 148 (2020), no. 7, 2783–2798.

\bibitem {NS2} Nathan Reading and Salvatore Stella. \textit{An affine almost positive roots model}. J. Comb. Algebra 4 (2020), no. 1, 1–59.
 
\bibitem {N} Nicolas Bourbaki. \textit{Lie groups and Lie algebras}. Chapters 4–6.  Elements of Mathematics (Berlin). Springer-Verlag, Berlin, 2002. xii+300 pp. ISBN: 3-540-42650-7.

\bibitem {SZ} Sergey Fomin and Andrei Zelevinsky, \textit{Cluster algebras. IV. Coefficients}. Compos. Math. 143 (2007), no. 1, 112–164.
\end{thebibliography}

\end{document}